\documentclass{amsart}
\usepackage[utf8]{inputenc}
\usepackage{amssymb,amsmath,mathtools,latexsym}
\usepackage{amsthm}
\usepackage{setspace}
\usepackage{graphicx}
\usepackage{mathrsfs}
\usepackage{mdwlist}
\usepackage{color}
\usepackage{bbm}
\usepackage{url}
\usepackage[colorlinks=true]{hyperref}
\usepackage{breakurl}
\usepackage[all,cmtip]{xy}
\usepackage[shortlabels]{enumitem}
\usepackage{tikz-cd}

\newtheorem{theorem}{Theorem}

\newtheorem{definition}[theorem]{Definition}

\newtheorem{conj}[theorem]{Conjecture}

\theoremstyle{remark}\newtheorem{remark}[theorem]{Remark}
\newtheorem{example}[theorem]{Example}
\newtheorem*{claim*}{Claim}

\numberwithin{theorem}{section}
\numberwithin{equation}{section}

\newcommand{\G}{\mathbf{G}} 

\newcommand{\SL}{\operatorname{SL}}
\newcommand{\SO}{\operatorname{SO}}
\newcommand{\Spin}{\operatorname{Spin}}
\newcommand{\torus}{\mathbf{T}} 

\newcommand{\A}{\mathbb{A}}

\newcommand{\disc}{\operatorname{disc}}
\renewcommand{\mod}{\, \mathrm{mod}\,} 
\newcommand{\norm}[1]{\|#1\|} 
\newcommand{\Ad}{\mathrm{Ad}}

\newcommand{\Bow}{\mathrm{Bow}}

\newcommand\rquot[2]{
  \mathchoice
  {
    \text{\raise0.5ex\hbox{$#1$}\big/\lower0.5ex\hbox{$#2$}}%
  }
  {
    #1\,/\,#2
  }
  {
    #1\,/\,#2
  }
  {
    #1\,/\,#2
  }
}

\newcommand\lquot[2]{
  \mathchoice
  {
    \text{\lower0.5ex\hbox{$#1$}\big\backslash\raise0.5ex\hbox{$#2$}}%
  }
  {
    #1\,\backslash\,#2
  }
  {
    #1\,\backslash\,#2
  }
  {
    #1\,\backslash\,#2
  }
}

\newcommand\lrquot[3]{
  \mathchoice
  {
    \text{\lower0.5ex\hbox{$#1$}\big\backslash\raise0.5ex\hbox{$#2$}\big/
      \lower0.5ex\hbox{$#3$}}%
  }
  {
    #1\,\backslash\,#2\,/\,#3
  }
  {
    #1\,\backslash\,#2\,/\,#3
  }
  {
    #1\,\backslash\,#2\,/\,#3
  }
}

\newcommand{\N}{\mathbb{N}}

\newcommand{\Q}{\mathbb{Q}}

\newcommand{\R}{\mathbb{R}}

\newcommand{\Z}{\mathbb{Z}}

\newcommand{\Hbf}{\mathbf{H}}

\newcommand{\id}{\mathrm{id}}

\DeclarePairedDelimiter\abs{\lvert}{\rvert}%

\makeatletter
\let\oldabs\abs
\def\abs{\@ifstar{\oldabs}{\oldabs*}}

\title{Representations of binary forms by quaternary quadratic forms}
\author{Wooyeon Kim, Andreas Wieser, Pengyu Yang}
\date{\today}
\thanks{
P.Y. is supported by National Key R\&D Program of China 2022YFA1007500.
W.K.~is supported by Korea Institute for Advanced Study, grant no.~HP101301.
A.W.~acknowledges the support of the Swiss National Science Foundation, grant no.~217944.
This material is based upon work supported by a grant from the Institute for Advanced Study School of Mathematics.
}

\begin{document}

\maketitle
\begin{abstract}
We prove a local-global principle for representations of binary by quaternary quadratic forms.
One of the main ingredients is a recent measure rigidity result of Einsiedler and Lindenstrauss for diagonalizable actions on quotients of products of $\mathrm{SL}_2$'s.
Based on this, it suffices to show that limits of the uniform measures on the associated rank one adelic toral packets have more entropy than one half of the maximal entropy.
The latter is proved using the Siegel mass formula and the determinant method as developed by Bombieri and Pila as well as Heath-Brown.
\end{abstract}

\section{Introduction}

Consider two integral quadratic forms $q:\Z^m \to \Z$ and $Q: \Z^n \to \Z$ for $m < n$.
We write $\disc(Q)$ for the discriminant of $Q$ i.e.~the determinant of any matrix representation of $Q$; $\disc(q)$ is defined analogously.
A \emph{representation} of $q$ by $Q$ is a linear map $\iota: \Z^m \to \Z^n$ with $Q\circ \iota = q$.
Furthermore, the representation $\iota$ is \emph{primitive} if $\iota(\Z^m) = (\iota(\Z^m)\otimes\Q) \cap \Z^n$.
Lastly, we will say that $q$ is (primitively) represented by $Q$ if a (primitive) representation of $q$ by $Q$ exists.

Clearly, a (primitive) representation can only exist if a local (primitive) representation $\iota_p: \Z_p^m \to \Z_p^n$ exists for every prime $p$. In short, we will say that $q$ is locally (primitively) represented by $Q$ if this local condition holds.

The \emph{local-global principle} asks whether $q$ is (primitively) represented when it is locally (primitively) represented.
This is in analogy to the classical Hasse-Minkowksi theorem which asserts that there is an isometry $(\Q^m,q) \to (\Q^n,Q)$ whenever there is an isometry over $\Q_p$ for every prime $p$ and over $\R$.
The Hasse-Minkowski theorem implies that $q$ is (primitively) represented by an element of the \emph{genus} of $Q$ if it is locally (primitively) represented by $Q$.
(Here, the genus of $Q$ is the set of integral quadratic forms locally equivalent to $Q$.)
The representation by an element of the genus of $Q$ may be upgraded to an element of the \emph{spin genus} of $Q$ whenever $n-m\geq 3$.
For indefinite forms the spin genus is trivial by work of Eichler \cite{eichler}, and hence so is the above local-global principle when $n-m\geq 3$.
We will from now on assume that $Q$ is positive definite.

\subsection{Representations of integers}
The local-global principle is particularly classical when $m=1$ in which case we are looking for (primitive) representations of numbers.
For $n\geq 4$, Kloosterman \cite{Kloosterman-diagonalforms} and Tartakovskii \cite{Tartakowskii1929} proved in the 1920's using the circle method that any sufficiently large number $D$ is primitively represented by $Q$ if it is locally primitively represented.
For $n=3$ the analogous claim is false, since local representation does not guarantee representation by the spin genus.
For instance, one can show elementarily that the quadratic form $x^2+xy+y^2+9z^2$ does not primitively represent any number of the form $4m^2$ where $m>0$ and $m\equiv 1 \mod 3$ (but it does so locally) --- see \cite[p.~115]{Watson-book}.
It is a general phenomenon that outside of finitely many square classes local representation implies representation by the spin genus.

Based on estimates for Fourier coefficients of half-integral weight forms by Duke \cite{duke88} and Iwaniec \cite{Iwaniec-halfintegral}, Duke and Schulze-Pillot \cite{DukeSchulzePillot} proved the following local-global principle:

\begin{theorem}[{\cite{DukeSchulzePillot}}]\label{thm:DukeSchulzePillot}
Suppose $n=3$.
If $D>0$ is sufficiently large and primitively represented by the spin genus of $Q$, then $D$ is primitively represented by $Q$.
\end{theorem}

A similar result has been established by Cogdell, Piatetski-Shapiro, and Sarnak \cite{Cogdell-squares} over number fields.
We also note that a weaker version of the above theorem can be obtained from the method of Linnik \cite{linnik} (see also \cite{ELMV-Ens,W-Linnik}); the ideas in the current article build on this method.

Lastly, we remark that none of the above methods for $n=3$ provide effective thresholds $D_0$ so that Theorem~\ref{thm:DukeSchulzePillot} holds for $D\geq D_0$ (due to potential Landau-Siegel zeroes).
This makes Theorem~\ref{thm:DukeSchulzePillot} difficult to use in explicit instances, see for example \cite{OnoSoundararajan}.
In contrast, effective thresholds do exist when $n \geq 4$, see e.g.~\cite{Hanke-DukeMathJ,SchulzePillot-explicitTartakovski}.

\subsection{Representations of quadratic forms}
Suppose now that $m>1$ i.e.~$q$ has at least two variables.
In the 70's, Hsia, Kitaoka, and Kneser \cite{HsiaKitaokaKneser} proved a local-global principle when $n\geq 2m+3$ under the additional condition that
\begin{align*}
\min(q) := \min_{x \in \Z^m\setminus\{0\}} q(x)
\end{align*}
is sufficiently large.
In a dramatic breakthrough 30 years later, Ellenberg and Venkatesh \cite{localglobalEV} improved the codimension assumption to $n-m \geq 5$ and under a splitting condition to $n-m \geq 3$.
Their work uses measure classification results from unipotent dynamics in an essential way --- specifically a $p$-adic variant of work of Mozes and Shah \cite{Mozes-Shah} based on measure classification results of Ratner~\cite{Ratner-measure,Ratner-padic} and Margulis, Tomanov \cite{MargulisTomanov} (see also Gorodnik, Oh \cite{GorOh-RationalPts}).
Recent work of Einsiedler, Lindenstrauss, Mohammadi, and the second named author \cite{effectivesemisimple} established effective equidistribution rates for semisimple adelic periods.
This work yields the following local-global principle when $n-m\geq 3$: 
when $q$ is locally primitively represented by $Q$ and $\min(q) \geq C \disc(Q)^A$ for effective constants $A,C>0$ depending only on $n$, then $q$ is primitively represented by $Q$.

When $n-m =2$, the above results from unipotent dynamics become inapplicable.
Nevertheless, the following is conjectured, also in analogy to the work of Duke and Schulze-Pillot \cite{DukeSchulzePillot} (Theorem~\ref{thm:DukeSchulzePillot} above).

\begin{conj}\label{conj:n-m=2}
Suppose $n-m =2$. If $q$ is primitively represented by the spin genus of $Q$ and $\min(q)$ is large enough (depending on $Q$), then  $q$ is primitively represented by $Q$.
\end{conj}

The current article makes progress on this conjecture in dimensions $m=2$ and $n=4$ i.e.~when $q$ is a binary integral quadratic form and $Q$ is a quaternary integral quadratic form (assumed to be positive-definite, as above).
In this case, we establish the following:

\begin{theorem}\label{thm:main}
Let $p_1,p_2$ be two distinct odd primes.
Then there exists $C=C(p_1,p_2,Q)>1$ with the following property.

Let $q$ be a primitive binary integral quadratic form such that $-\disc(q)\disc(Q)$ is a non-zero square modulo $p_1,p_2$.
If $q$ is primitively represented by the spin genus of $Q$ and $\min(q) \geq C$ then $q$ is primitively represented by $Q$.
\end{theorem}

We note that the two auxiliary primes $p_1,p_2$ are an artifact of our method. Indeed, they make the problem accessible by known classification results for measures invariant and ergodic under higher-rank diagonalizable actions. Specifically, we use recent work of Einsiedler and Lindenstrauss \cite{EL-nonmaximal} for irreducible quotients of a product of $\SL_2$'s.
The additional assumption that $q$ be primitive can be weakened and should not be seen as central to our approach.

Our methods also yield lower bounds on the number of primitive representations $r(q,Q)$ of $q$ by $Q$.
Denote by $r(q,\mathrm{spin}(Q))$ the usual weighted average of the primitive representation numbers over the spin genus.
We establish the following stronger version of the above theorem:

\begin{theorem}\label{thm:main2}
There exists $\delta>0$ with the following property.
Let $p_1,p_2$ be two distinct primes and let $q,Q$ be as in Theorem~\ref{thm:main}.
Then
\begin{align*}
r(q,Q) \geq r(q,\mathrm{spin}(Q)) \big( \delta + \varepsilon(q))
\end{align*}
where $\varepsilon(q)$ is a function depending on $Q,p_1,p_2$ and on $q$ with $\varepsilon(q) \to 0$ as $\min(q)$ goes to infinity.
\end{theorem}

As already alluded to, the technology available to attack Conjecture~\ref{conj:n-m=2} is generally much weaker than in the works mentioned earlier for the local-global principle when $n-m \geq 3$.
For instance, the effective equidistribution statements proven in \cite{effectivesemisimple} allow for a precise asymptotic statement for representation numbers $r(q,Q)$ with an effective estimate on $\varepsilon(q)$ in terms of $\min(q)$. Statements of this kind are very far out of reach in the current setting where we do not establish equidistribution results.

We also remark that the local-global principle (in the form of Conjecture~\ref{conj:n-m=2}) for representations of binary by quaternary quadratic forms is very closely related to the mixing conjecture of Michel and Venkatesh \cite{MVIHES}.
The mixing conjecture has seen striking progress in recent years with works of Khayutin \cite{IlyaAnnals}, Blomer and Brumley~\cite{BlomerBrumley-mixing}, Blomer, Brumley, and Khayutin \cite{BlomerBrumleyKhayutin}, as well as Blomer, Brumley, Radziwi\l\l.
The differing methods in these works rely on underlying product structure that is unavailable for the problem of the present article.
Correspondingly, there are currently no viable analytic approaches to the local-global principle for $m=2,n=4$.
An interesting result based on analytic methods is due to Schulze-Pillot \cite{SchulzePillot-2in4}, who shows that a positive proportion of binary quadratic forms of discriminant $\disc(q)$ are represented.

In the remainder of this announcement, we will present ideas of the proof for Theorem~\ref{thm:main}. As mentioned, our approach is influenced by Linnik's method \cite{linnik} as reinterpreted in \cite{ELMV-Ens,W-Linnik} and we will point out similarities and differences.

\medskip
\textbf{Acknowledgments:}
The authors are grateful towards Manfred Einsiedler, Elon Lindenstrauss, Peter Sarnak, and Akshay Venkatesh for their interest in this project and for their encouragement.
We also thank Menny Aka, Farrell Brumley, Zhizhong Huang, Ilya Khayutin, Rainer Schulze-Pillot, Per Salberger, Ye Tian, and Katherine Woo for fruitful discussion on various topics.
Last but not least, we thank the Forschungsinstitute f\"ur Mathematik at ETH Zurich, the Institute for Advanced Study, the Korea Institute for Advanced Study, the Morningside Center of Mathematics at CAS, and the Simons-Laufer Institute for providing an excellent work environment.

\section{Ideas of the proof}

\subsection{Density of toral periods}

We begin by reinterpreting our goal in terms of certain adelic periods.
Denote by $\G = \Spin_Q$ the spin group of the quadratic form $Q$.
By definition, $\G$ is the universal cover of $\SO_Q$ and we denote by $\rho: \G\to\SO_Q$ the covering map.
Note that $\G$ is a $\Q$-anisotropic $\Q$-form of $\SL_2\times \SL_2$.
One can show that $\G$ is an inner form if and only if the discriminant $\disc(Q)$ is a rational square; this should be seen as the `atypical' case (also in view of earlier comments on the mixing conjecture).
Define the compact space
\begin{align}\label{eq:localglobal-ambientspaces}
X = \SO_Q(\Q)\rho(\G(\A)) \subset \lquot{\SO_Q(\Q)}{\SO_Q(\A)}
\end{align}
and let $\mu_X$ be the unique $\rho(\G(\A))$-invariant probability measure on $X$.

Suppose now that we are given a sequence of binary quadratic forms $q_i$ with $\min(q_i) \to \infty$ such that $-\disc(q_i)\disc(Q)$ is a non-zero square modulo $p_1,p_2$.
We may assume that $Q$ primitively represents each $q_i$ and show instead that, for large enough $i$, any element of the spin genus of $Q$ also primitively represents $q_i$.

For each $i$, let $\iota_i$ be a primitive representation of $q_i$ by $Q$ and define the $\Q$-torus
\begin{align*}
\torus_i &= \{g \in \G: g.\iota_i(x) = \iota_i(x) \text{ for all }x \in \Z^2\}.
\end{align*}
These are one-dimensional $\Q$-anisotropic tori.

\begin{remark}
If $\G$ is an inner form of $\SL_2 \times \SL_2$, we may write $\G = \mathbf{B}^1 \times \mathbf{B}^1$ where $\mathbf{B}$ is a quaternion algebra over $\Q$ ramified at the infinite place and where $\mathbf{B}^1$ is the group of norm one units. In this case, one may verify that $\torus_i$ projects non-trivially to both factors of $\G$.
\end{remark}

Define the adelic toral periods
\begin{align*}
Y_i = \SO_Q(\Q) \rho(\torus_i(\A)) \subset X
\end{align*}
and let $\nu_i$ be the $\rho(\torus_i(\A))$-invariant probability measure on $Y_i$.

\begin{theorem}\label{thm:periods}
There exists $\delta>0$ absolute with the following property.
Let $\nu$ be any weak${}^\ast$-limit of the measures $\nu_i$.
Then $\nu \geq \delta \mu_X$.
\end{theorem}

In particular, this shows that the sequence $Y_i$ is asymptotically dense i.e.~given any open set $\mathcal{O}\subset X$ there exists $i_0$ so that $Y_i \cap \mathcal{O} \neq\emptyset$ for all $i\geq i_0$.
The deduction of our main theorems from Theorem~\ref{thm:periods} is fairly standard and appears in a similar form already in \cite{localglobalEV,effectivesemisimple}.
We will therefore focus on outlining the proof of Theorem~\ref{thm:periods} from now on.

When $\min(q_i)$ grows much slower than the discriminants $\disc(q_i)$, one can apply an argument known as `equidistribution in stages', among other names.
In this situation, one obtains equidistribution of the measures $\nu_i$, which implies, for instance, exact asymptotics in Theorem~\ref{thm:main2}.
To briefly outline the argument for equidistribution, suppose $v_i \in \iota_i(\Z^2)$ is a shortest vector (with $Q$-value $\min(q_i)$) and let $\Hbf_i = \{g\in \G:g.v_i = v_i\}$. Then the adelic period $Z_i = \SO_Q(\Q) \rho(\Hbf_i(\A))$ has volume polynomial in $\min(q_i)$ and one may apply equidistribution of $Y_i$ in $Z_i$ and of $Z_i$ in $X$.
The former is a version of Duke's equidistribution theorem \cite{duke88} and its variants, see e.g.~\cite[\S4]{ELMVAnn} and the references therein.
The latter follows from recent progress on effective equidistribution of semisimple adelic periods; see \cite{EMMV,effectivesemisimple}.
(Note that if $\G$ is an inner form of $\SL_2\times \SL_2$ this is a version of effective equidistribution of Hecke points proven earlier in \cite{ClozelOhUllmo}).

As the above argument is relatively standard, we will assume from now on that
\begin{align}\label{eq:min not small}
\min(q_i) \geq |\disc(q_i)|^{\eta}
\end{align}
for some $\eta>0$.

\subsection{Applying measure rigidity results}
Recall our splitting assumptions at the primes $p_1,p_2$ in Theorem~\ref{thm:main}.
These imply that the tori $\torus_i$ split over $\Q_{p_1}$ and $\Q_{p_2}$.
To unify the invariance at $p_1,p_2$, we may pick a bounded sequence $k_i \in \G(\Q_{p_1})\times \G(\Q_{p_2})$ such that $A = k_i^{-1} (\torus_i(\Q_{p_1})\times \torus_i(\Q_{p_2}))k_i$ does not depend on $i$.
Set $Y_i' = Y_i k_i$ and let $\nu_i'$ be the Haar probability measure on $Y_i'$.
It is sufficient to show that any weak${}^\ast$-limit $\nu'$ of the measures $\nu_i'$ satisfies the conclusion of Theorem~\ref{thm:periods}.
Clearly, $\nu'$ is $A$-invariant.

In a series of fundamental works including \cite{EL-maximal,EL-lowentropy,EL-symmetry,EL-joiningsPIHES}, Einsiedler and Lindenstrauss have successfully classified, in many instances, probability measures on homogeneous spaces invariant and ergodic under higher-rank diagonalizable actions, assuming positive entropy.
Most recently, they resolved in \cite{EL-nonmaximal} a first instance of this broad program when the action is non-maximal, specifically for actions on irreducible quotients of products of $\SL_2$'s.
Suppose for simplicity that $\G$ is an outer form of $\SL_2\times \SL_2$ (as always compact over $\R$).
By~\cite{EL-nonmaximal}, for any ergodic component $\mu$ of $\nu'$ we have that either
\begin{itemize}
    \item $\mu$ is homogeneous or
    \item $h_\mu(a) =0 $ for any $a \in A$.
\end{itemize}
Here, $h_\mu(a)$ denotes the entropy of $\mu$ with respect to $a$.
The homogeneous measures are easily analyzed in our setting and one can verify that any homogeneous measure $\mu$ with positive entropy satisfies either $h_\mu(a) = \frac{1}{2}h_{\mu_X}(a)$ (in which case $\mu$ is the homogeneous measure for an intermediate period) or $\mu = \mu_X$.
In summary, the above application of \cite{EL-nonmaximal} yields the following 

\textbf{Conclusion:}
It suffices to show that for some $a \in A$ we have
\begin{align}\label{eq:entropy > 0.5}
h_{\nu'}(a) > \tfrac{1}{2}h_{\mu_X}(a).
\end{align}

We establish \eqref{eq:entropy > 0.5} by an argument similar to the aforementioned approach of Linnik e.g.~equidistribution of integer points on spheres.
Fix $a \in A$.
Choose an open neighborhood $B \subset \G(\A)$ of the identity and define, for any $n\geq 0$, the two-sided Bowen balls
\begin{align*}
\Bow_n = \bigcap_{|k| \leq n} a^k B a^{-k}.
\end{align*}
If there is $\delta>0$ and a sequence $n_i \to \infty$ for which the self-correlation estimate
\begin{align}\label{eq:selfcorrelation}
\nu_i' \times \nu_i'\big(\{ (x,y) \in X^2: y \in x \rho(\Bow_{n_i})\}\big)
\leq \mathrm{e}^{-2n_i(\frac{1}{2}+\delta) h_{\mu_X}(a)}
\end{align}
holds, then $h_{\nu'}(a) \geq (\frac{1}{2}+\delta)h_{\mu_X}(a)$ and \eqref{eq:entropy > 0.5} follows.

\subsection{A counting problem}
It remains to prove a self-correlation estimate of the shape in \eqref{eq:selfcorrelation}.
An estimate of this type typically corresponds to a counting problem.
As \eqref{eq:selfcorrelation} is a statement about an individual period, we write $q=q_i$ for simplicity. Set $p:=p_1$.
With a suitable choice of $a\in A$, the counting problem in our setting is given as follows.

Write $q(x,y) =Ax^2+Bxy+Cy^2$ and set $D=\disc(q) = AC-\frac{1}{4}B^2$.
For any $n \geq 1$, consider the set $\mathcal{X}(n)$ of pairs $(\iota_1,\iota_2)$ with the following properties:
\begin{itemize}
\item $\iota_1,\iota_2$ are primitive representations of $q$ by $Q$ with distinct images.
    \item There exists an rotation $k$ in the plane $\iota_1(\Q_p^2)$ (that is, $k \in \SO_{Q}(\Z_p)$ with $k|_{\iota_1(\Q^2)^\perp}=\id$) such that
    \begin{align*}
    \iota_2(x) \equiv k\iota_1(x) (\mod p^{2n})
    \end{align*}
    for all $x\in \Z^2$.
\end{itemize}
Then \eqref{eq:selfcorrelation} amounts to showing for some $n$ with $p^{(2+\delta)n} \leq D^{1/4}$ (and with $n\to \infty$ as $D \to \infty$)
\begin{align}\label{eq:bound on pairs}
\#\mathcal{X}(n) \ll_\varepsilon \frac{D^{1+\varepsilon}}{p^{(4+2\delta)n}}.
\end{align}
We note that, more precisely, the left-hand side of \eqref{eq:selfcorrelation} for time $n$ is bounded by $\ll_\varepsilon D^\varepsilon (\frac{1}{\sqrt{D}}+ \frac{1}{D}\#\mathcal{X}(n))$ where the first summand is the `diagonal' contribution.
For the above implicit choice of $a\in A$, the $\Ad$-eigenvalues are $p^{\pm 2}$ and, correspondingly, the maximal entropy is $h_{\mu_X}(a) = 4\log(p)$.

\begin{remark}
For Linnik-type problems such as Theorem~\ref{thm:DukeSchulzePillot} with a splitting condition at a prime $p$, one counts, given a positive definite ternary integral quadratic form $Q$ and $D>0$, the number of pairs $(v,w)$ where $v,w \in \Z^3$ are primitive with $Q(v)=Q(w)=D$ and $v \equiv w \mod p^n$. See for example \cite{ELMV-Ens,W-Linnik}, where good estimates are obtained using the Siegel mass formula.
\end{remark}

We now discuss how to obtain an estimate as in \eqref{eq:bound on pairs}.
A crucial application of the Siegel mass formula \cite{Siegel-I,Siegel-II,Siegel-III} implies that it suffices to count all possible quadratic forms attainable on the sublattice $\iota_1(\Z^2)+\iota_2(\Z^2)$.
For the purposes of this outline, we focus on the `generic' case where $\iota_1(\Z^2)+\iota_2(\Z^2)$ is a rank $4$ lattice.
Define the half-integers
\begin{equation*}
\begin{split}
x_1 = \langle \iota_1(e_1),\iota_2(e_1)\rangle_{Q},\ x_2 = \langle \iota_1(e_1),\iota_2(e_2)\rangle_{Q},\\
x_3 = \langle \iota_1(e_2),\iota_2(e_1)\rangle_{Q},\ x_4 = \langle \iota_1(e_2),\iota_2(e_2)\rangle_{Q}.
\end{split}
\end{equation*}
Thus, the quadratic form on $\iota_1(\Z^2)+\iota_2(\Z^2) = \Z\iota_1(e_1)+\Z\iota_1(e_2)+\Z\iota_2(e_1)+\Z\iota_2(e_2)$ is given by
\begin{align*}
\begin{pmatrix}
A & B/2 & x_1 & x_2 \\
B/2 & C & x_3 & x_4 \\
x_1 & x_3 & A & B/2 \\
x_2 & x_4 & B/2 & C
\end{pmatrix}.
\end{align*}
Note that the determinant of the above symmetric matrix is the determinant of $Q$ times the square of the index $x_0 = [\Z^4:\iota_1(\Z^2)+\iota_2(\Z^2)]$.
In other words, the tuple $(x_0,x_1,\ldots,x_4)$ is a point on the variety
\begin{equation}\label{eq:degree4 equation}
\begin{split}
\disc(Q)x_0^2 &= (x_1x_4-x_2x_3 -D)^2 - (Cx_1-Ax_4)^2 \\
&\qquad- (Bx_1-A(x_2+x_3))(Bx_4-C(x_2+x_3)).
\end{split}
\end{equation}
The tuple $(x_0,x_1,\ldots,x_4)$ satisfies further restrictions and we summarize all the information in the following definition.

\begin{definition}\label{def:set S}
For any $n\in \N$ the subset $\mathcal{S}(n) \subset \frac{1}{2}\Z^5$ is the set of tuples $x =(x_0,\ldots,x_4)$ satisfying \eqref{eq:degree4 equation} as well as the following constraints:
\begin{itemize}
\item $|x_1| \leq A$, $|x_2|,|x_3| \leq \sqrt{AC}$, and $|x_4| \leq C$.
\item We have
\begin{equation}\label{eq:originalcongeqforx}
    \begin{aligned}
x_1x_4-x_2x_3 &\equiv D \mod p^{4n},\\
Cx_1-Ax_4 &\equiv 0 \mod p^{4n},\\
C(x_2+x_3)-Bx_4 &\equiv 0 \mod p^{4n},\\
A(x_2+x_3)-Bx_1 &\equiv 0 \mod p^{4n}.
\end{aligned}
\end{equation}

\item If $x_1 =A$ or $x_4 = C$ then $x_0 = 0$, $x_2=x_3=\frac{B}{2}$.
If $x_1 =-A$ or $x_4 = -C$ then $x_0 = 0$, $x_2=x_3=-\frac{B}{2}$.
\end{itemize}
\end{definition}

Here, the first and third bullet points are a simple consequence of the Cauchy-Schwarz inequality while \eqref{eq:originalcongeqforx} is a consequence of the congruence condition in the definition of $\mathcal{X}(n)$.

A `trivial' bound on the size of $\mathcal{S}(n)$ is given by roughly $\frac{D^2}{p^{12n}}$. This amounts to having at least one half of the maximal entropy (i.e.~non-strict inequality in \eqref{eq:entropy > 0.5}) and does not invoke \eqref{eq:degree4 equation} meaningfully.
It remains to improve on this trivial bound.
To that end, we use the determinant method developed by Bombieri and Pila \cite{BombieriPila}, Heath-Brown \cite{HeathBrown-annals}, and Salberger \cite{Salberger-PLMS}.
Specifically, we use a variant of the result of Bombieri and Pila \cite{BombieriPila} counting integral points on a rational planar curve; the variant can be established using Heath-Brown's $p$-adic approach from \cite{HeathBrown-annals}.
In the following we outline our argument in more detail.

We choose $p^n$ to be close to $D^{\frac{1}{8}}$, though slightly smaller by a gap in the exponent.
With this choice, we need to show that 
\begin{align*}
\#\mathcal{S}(n) \ll \sqrt{D}.
\end{align*}
We note that to leverage \eqref{eq:degree4 equation} geometric information on the variety $\mathbf{V}$ cut out by \eqref{eq:degree4 equation} has to be used.
For instance, affine subspaces of low height contained in $\mathbf{V}$ can be potentially problematic, or more generally low degree subvarieties.
This is reflected in our application of the determinant method, where the following two examples take a special role.

\begin{example}
$\mathbf{V}$ contains the affine linear subvariety cut out by $x_1=A$, $x_2=x_3=\frac{B}{2}$, and $x_0=0$, on which there are $\gg \frac{C}{p^{4n}}$ many points of $\mathcal{S}(n)$.
In view of \eqref{eq:min not small}, this subset of $\mathcal{S}(n)$ is not problematic.
\end{example}

\begin{example}
If $\disc(Q) = d_0^2$ for some $d_0$, the variety
$\mathbf{V}$ contains the subvariety cut out by $Cx_1=Ax_4$, $Bx_1= A(x_2+x_3)$, and $d_0 x_0 = x_1x_4-x_2x_3 -D$. Our assumption that $q$ be primitive implies that $(x_1,x_2+x_3,x_4)$ is a multiple of $(A,B,C)$ and, in particular, the subvariety does not contribute meaningfully.

For the opposite extreme, if for instance $q(x,y) = A(x^2+y^2)$ then $x_1=x_4$ and $x_2+x_3=0$. 
There are $\asymp A$ choices for $x_1$ and, by the quadratic congruence condition in \eqref{eq:originalcongeqforx}, around $Ap^{-4n}$ choices for $x_2$ so that the total number of points of $\mathcal{S}(n)$ on this subvariety is $\asymp A^2p^{-4n} =Dp^{-4n}$. 
This is strictly larger than our desired bound for $\mathcal{S}(n)$ and illustrates the role of our primitivity assumption on $q$.
\end{example}

To restrict the following argument to the most interesting case, we assume henceforth that the quadratic form $q$ is `balanced' i.e.~$\min(q) \gg \sqrt{D}$. We also assume that $q$ is reduced so that $A = \min(q)$; this latter assumption does not restrict the generality.

We begin by linearizing the quadratic equation in \eqref{eq:originalcongeqforx}.
Thus, we fix a solution $w' = (w_1,\ldots,w_4)$ of \eqref{eq:originalcongeqforx} modulo $p^{2n}$ (there are $\ll p^{2n}$ many such classes) and wish to count all the points $x \in \mathcal{S}(n)$ for which $x' = (x_1,\ldots,x_4)$ reduces to $w$ modulo $p^{2n}$.
The congruence condition \eqref{eq:originalcongeqforx} in the new coordinates $(y_1,
\ldots,y_4)$ given by $x_i = w_i + p^{2n}y_i$ translate to
\begin{equation}\label{eq:congeqfory}
    \begin{aligned}
w_4y_1+w_1y_4-w_2y_3-w_3y_2 &\equiv 0 \mod p^{2n},\\
Cy_1-Ay_4 &\equiv 0 \mod p^{2n},\\
C(y_2+y_3)-By_4 &\equiv 0 \mod p^{2n},\\
A(y_2+y_3)-By_1 &\equiv 0 \mod p^{2n}.
\end{aligned}
\end{equation}
We may assume that $w'$ lifts to (the last four coordinates of) an element of $\mathcal{S}(n)$ as otherwise $\mathcal{S}(n)$ contains no points in the fiber above $w'$.
In particular, we have
\begin{align*}
|y_1|,|y_2|,|y_3|,|y_4| \ll \sqrt{D}p^{-2n}.
\end{align*}

The congruence equations in \eqref{eq:congeqfory} are the desired linearized congruence equations. They cut out a sublattice $\Lambda_w$ of $\frac{1}{2}\Z^4$ of index expected to be around $p^{6n}$. More precisely, the index is at least $p^{6n-\nu}$ if $p^{\nu+1}\nmid (w_2-w_3)$.
The set of $w$ with $p^{\nu+1}\mid (w_2-w_3)$ is a smaller set and can be removed from consideration. For simplicity, we suppose here $p \nmid (w_2-w_3)$.

Let $v_1,\ldots,v_4$ denote a Minkowski reduced basis of $\Lambda_w$ with respect to the usual Euclidean norm. 
In particular, $\norm{v_1} \cdots \norm{v_4} \asymp p^{6n}$ and $\norm{v_1}\leq \norm{v_2}\leq\norm{v_3}\leq\norm{v_4}$.
We can write
\begin{align*}
y' = (y_1,\ldots,y_4) = z_1 v_1+z_2 v_2+z_3 v_3+z_4 v_4
\end{align*}
for new variables $z_1,\ldots,z_4$ with $|z_i| \leq B_i$ where $B_i \asymp \sqrt{D}p^{-2n}\norm{v_i}^{-1}$. 
Observe that $B_i \gg 1$ as $\norm{v_4} \leq p^{2n}$.

Generically, our argument now proceeds as follows.
Insert the above new coordinates $(z_1,\ldots,z_4)$ into \eqref{eq:degree4 equation} and fix $z_2,z_3,z_4$. The so-obtained equation ought to be an irreducible planar curve and the number of points on it is bounded by $B_1^{1/2}$ using the aforementioned variant of the result of Bombieri-Pila in \cite{BombieriPila}.
This yields a saving of $B_1^{1/2} \geq D^{\frac{1}{16}}$ which suffices for our purposes.
We observe the following problems with this argument:
\begin{itemize}
\item The above equation in $x_0,z_4$ might not be irreducible or, equivalently, the right-hand side of \eqref{eq:degree4 equation} might be a square as a function in $z_4$. One shows that this can only happen for few values of $z_2,z_3,z_4$.
\item If $B_2 \leq D^{\delta}$ for some very small $\delta>0$, then we only consider few values of $z_1,z_2,z_3$ and the information of the previous bullet point is not useful. In this case where the lattice $\Lambda_w$ has a very short vector $v_1$ we invoke \eqref{eq:congeqfory} again to show that $w'$ is significantly restricted, and obtain a gain in the count of these `bad' $w'$.
\end{itemize}
Overall, the above analysis completes our outline showing $\#\mathcal{S}(n) \ll \sqrt{D}$. As explained, this establishes more than one half of maximal entropy and, with it, Theorem~\ref{thm:periods}.

\bibliographystyle{plain}
\bibliography{Bibliography}

\end{document}